\documentclass[letterpaper,english,aps,prl,superscriptaddress,twocolumn,longbibliography]{revtex4-1}
\usepackage[T1]{fontenc}
\usepackage[latin9]{inputenc}
\usepackage{array}
\usepackage{booktabs}
\usepackage{amsmath}
\usepackage{graphicx}

\makeatletter

\pdfpageheight\paperheight
\pdfpagewidth\paperwidth

\providecommand{\tabularnewline}{\\}

\usepackage{graphicx}
\usepackage{bm}
\usepackage{xcolor}
\usepackage{amsfonts}
\usepackage{amsmath}
\usepackage{amssymb}
\usepackage{url}
\usepackage[colorlinks,citecolor=blue,linkcolor=blue]{hyperref}

\makeatother

\usepackage{babel}
\begin{document}
\title{Geodesic bound of the minimum energy expense to achieve membrane separation
within finite time}
\author{Jin-Fu Chen}
\address{Beijing Computational Science Research Center, Beijing 100193, China}
\address{Graduate School of China Academy of Engineering Physics, No. 10 Xibeiwang
East Road, Haidian District, Beijing, 100193, China}
\address{now at: School of Physics, Peking University, Beijing 100871, China}
\author{Ruo-Xun Zhai}
\address{Graduate School of China Academy of Engineering Physics, No. 10 Xibeiwang
East Road, Haidian District, Beijing, 100193, China}
\author{C. P. Sun}
\address{Beijing Computational Science Research Center, Beijing 100193, China}
\address{Graduate School of China Academy of Engineering Physics, No. 10 Xibeiwang
East Road, Haidian District, Beijing, 100193, China}
\author{Hui Dong}
\email{hdong@gscaep.ac.cn}

\address{Graduate School of China Academy of Engineering Physics, No. 10 Xibeiwang
East Road, Haidian District, Beijing, 100193, China}
\date{\today}
\begin{abstract}
To accomplish a task within limited operation time typically requires
an excess expense of energy, whose minimum is of practical importance
for the optimal design in various applications, especially in the
industrial separation of mixtures for purification of components.
Technological progress has been made to achieve better purification
with lower energy expense, yet little is known about the fundamental
limit on the least excess energy expense in finite operation time.
We derive such a limit and show its proportionality to the square
of a geometric distance between the initial and final states and inverse
proportionality to the operation time $\tau$. Our result demonstrates
that optimizing the separation protocol is equivalent to finding the
geodesic curve in a geometric space. Interestingly, we show the optimal
control with the minimum energy expense is achieved by a symmetry-breaking
protocol, where the two membranes are moved toward each other with
different speeds.
\end{abstract}
\maketitle
Separating components in a chemical mixture into their purer forms
is critical in industrial applications \citep{Sholl2016}, such as
water purification \citep{Koros2012,Alvarez2018,Noamani2019}, pharmaceutical
and biological industry \citep{Reis2007,Xie2008}, and the environmental
science of $\mathrm{CO}_{2}$ controls \citep{FalkPedersen1997,Favre2007,Merkel2010,Adewole2013}.
The membrane separation is one of the most promising technologies
due to its low energy consumption and environmental-friendly operation
\citep{Britan1991,Koros1993,Spillman1995,Pabby2015,Castel2018,Purkait2018}.
Significant efforts have been made to improve the efficiency of the
membranes from the perspective of material design \citep{Pabby2015,Purkait2018}.
Nevertheless, the extent to which these properties benefit the energy
expense in separation processes remains elusive. An important question
arising naturally is whether there is a lower bound to the energy
consumption posted by the basic laws of thermodynamics. If so, such
a bound shall shed light on the material synthesis and the protocol
design of separation processes.

We seek a fundamental bound of the least energy expense in finite
operation time as a consequence of the basic law of thermodynamics,
particularly the geometric structure of the thermal equilibrium configuration
space \citep{Weinhold1975,Ruppeiner1979,Salamon1983,Crooks2007,Sivak2012,Scandi2018,Chen2021,Li2022}.
For quasi-static processes with infinite operation time, the basic
laws of thermodynamics have already provided a universal lower bound
for the energy consumption: the performed work should exceed the free
energy change of the separated final state and the mixed initial state
$W\geq\Delta F$ \citep{book:18204,book:18254}. For a practical finite-time
separation process, an excess work beyond the quasi-static separation
is typically required \citep{book:18254,Tsirlin2002,Castel2018} and
should be optimized to reach the fundamental limit \citep{Xu1996,Tsirlin2002,Sieniutycz2017,Sieniutycz2018}.
We convert such a task of finding the minimum energy expense into
finding the shortest path in a configuration space, demonstrating
the first application of geometric optimization in the membrane separation
process.

\begin{figure*}
\includegraphics{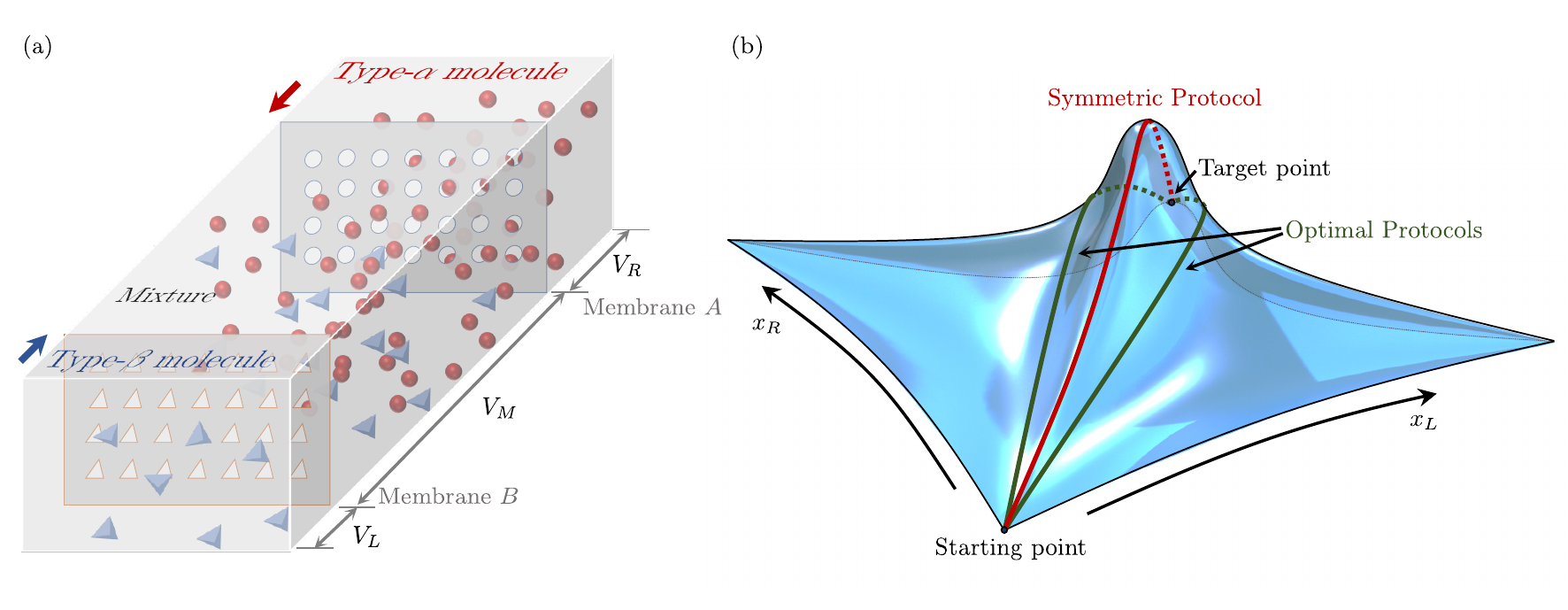}

\caption{Separation of binary mixed gases.\textbf{ }(a) The schematic of the
separation process. The mixture consists of two types of molecules,
type-$\alpha$ (red balls) and type-$\beta$ (blue triangular pyramids).
The semipermeable membrane A (B) is designed for only type-$\alpha$($\beta$)
to pass through. The chamber is divided into three compartments, the
left one (with volume $V_{L}$) for the purified type-$\alpha$ molecules,
the middle one (with volume $V_{M}$) for the molecular mixture, and
the right one (with volume $V_{R}$) for the purified type-$\beta$
molecules. The two membranes are pushed towards each other during
the separation process. (b) Sketch of the landscape of the configuration
space spanned by $x_{k}=V_{k}/V_{t}$ ($k=L,R$) for the separation
process. The metric is reflected by the surface in the three-dimensional
space with an example of the equally mixed gases with symmetric diffusion
coefficient. The excess work is reflected by the length of the path
between the start point $\left(0,0\right)$ and the target point $\left(0.5,0.5\right)$.}

\label{fig:scheme}
\end{figure*}

\textit{Modeling the membrane separation.} Consider a binary mixture
of molecules in a chamber with volume $V_{t}$, as illustrated in
Fig. \ref{fig:scheme}(a). The system is immersed in a thermal bath
with temperature $T_{0}$. Two species of the molecules are shown
as red balls (type-$\alpha$) and blue triangular pyramids (type-$\beta$)
with the numbers of the molecules $N_{\alpha}$ and $N_{\beta}$.
The separation process is performed by mechanically moving two semipermeable
membranes A and B from the two ends towards each other. The membrane
A (B) is designed with the properties allowing only type-$\alpha$
(type-$\beta$) molecules to penetrate. The chamber is divided into
three compartments with the two membranes: the left one with volume
$V_{L}$ for the purified type-$\alpha$ molecules, the middle one
with volume $V_{M}$ for the molecular mixture, and the right one
with volume $V_{R}$ for the purified type-$\beta$ molecules. The
three volumes satisfy the condition $V_{L}+V_{M}+V_{R}=V_{t}$. At
the end of the separation, the two membranes A and B contact each
other, i.e. $V_{M}=0$, and the two gases are purified in the left
and the right compartments.

We denote the number of type-$\sigma$ ($\sigma=\alpha,\beta$) molecules
in the compartment $k$ ($k=L,M,R$) as $N_{\sigma k}$. The impermeability
of membrane A (B) to type-$\alpha$ (type-$\beta$) molecules results
in $N_{\beta L}=N_{\alpha R}=0$. For the gaseous molecules, the relaxation
time is far shorter than the operation time, typically referenced
as the endo-reversible region \citep{Curzon1975,Salamon1981}, where
the status of the gas can be described by macroscopic parameters,
i.e., the pressure and the temperature. We assume the heat exchange
inside the system is fast enough to ensure a global temperature $T$
for the molecules in the whole chamber. It is worth noting that the
temperature $T$ typically differs from the bath temperature $T_{0}$
due to the finite-time operation. The equation of state for the gases
in each compartment is
\begin{equation}
p_{\sigma k}V_{k}=N_{\sigma k}k_{B}T,\label{eq: idealgas}
\end{equation}
where $k_{B}$ is the Boltzmann constant, and $p_{\sigma k}$ is the
corresponding partial pressure.

Two types of relaxations exist in the separation process, the particle
transport across the membranes and the heat conduction through the
wall of the chamber. According to Fick's law \citep{Hille2001}, the
particle transport flux $\mathcal{J}_{\alpha}$ ($\mathcal{J}_{\beta}$)
of type-$\alpha$ (type-$\beta$) molecules across the membrane A
(B) is proportional to the particle-density difference $\Delta c_{\alpha}=N_{\alpha M}/V_{M}-N_{\alpha L}/V_{L}$
($\Delta c_{\beta}=N_{\beta M}/V_{M}-N_{\beta R}/V_{R}$) across the
membrane A (B), namely, $\mathcal{J}_{\alpha}=\mu_{\alpha}\Delta c_{\alpha}$
and $\mathcal{J}_{\beta}=\mu_{\beta}\Delta c_{\beta}$, where $\mu_{\alpha}$
($\mu_{\beta}$) is the diffusion coefficient of the type-$\alpha$
(type-$\beta$) molecules penetrating membrane A (B). The change rates
of the molecular numbers are described by $\dot{N}_{\alpha L}=\mathcal{J}_{\alpha}\mathcal{A}$
and $\dot{N}_{\beta R}=\mathcal{J}_{\beta}\mathcal{A}$, and are explicitly
\begin{equation}
\begin{split}\dot{N}_{\alpha L}= & \mu_{\alpha}\mathcal{A}\left(\frac{N_{\alpha M}}{V_{M}}-\frac{N_{\alpha L}}{V_{L}}\right)\,,\\
\dot{N}_{\beta R}= & \mu_{\beta}\mathcal{A}\left(\frac{N_{\beta M}}{V_{M}}-\frac{N_{\beta R}}{V_{R}}\right)\,,
\end{split}
\label{eq:particle transport}
\end{equation}
where $\mathcal{A}$ is the area of membranes, assumed the same for
both membranes A and B. We adopt Newton's law of cooling to describe
the heat conduction between the system and the thermal bath, i.e.
$\dot{Q}=-C_{V}\gamma(T-T_{0})\,$. Here $C_{V}$ is the heat capacity
of the system at constant volume, e.g. $C_{V}=3/2N_{t}k_{B}$ for
the ideal single-atom gas. $N_{t}=N_{\alpha}+N_{\beta}$ is the total
number of the two types of molecules. $\gamma$ is the cooling rate
of the system. With the first law of thermodynamics, the evolution
of the gas temperature $T$ is described as
\begin{equation}
C_{V}\dot{T}=-C_{V}\gamma(T-T_{0})+\dot{W}\,,\label{eq:temperature}
\end{equation}
where $\dot{W}$ is the mechanical work rate performed on the gas
while moving the membranes.

\textit{Excess energy expense.} To accomplish the separation, the
mechanical work performed by moving the membranes is 
\begin{equation}
W=-\int_{0}^{\tau}\sum_{\sigma=\alpha,\beta}\sum_{k=L,M,R}p_{\sigma k}\dot{V}_{k}dt\,,\label{eq:W_accurate}
\end{equation}
where $\tau$ is the operation time of the separation process and
$\dot{V}_{k}=dV_{k}/dt$ is the volume change rate of the compartment
$k$. For convenience, we define two dimensionless parameters determining
the configuration of the system, $x_{L}\equiv V_{L}/V_{t}$ and $x_{R}\equiv V_{R}/V_{t}$.
The two species of molecules are completely mixed initially with $x_{L}(0)=x_{R}(0)=0$.
They are separated finally with $x_{L}(\tau)$ and $x_{R}(\tau)$
satisfying $x_{L}(\tau)+x_{R}(\tau)=1$. During the separation, there
is a constraint condition $0\leq x_{L}(t)+x_{R}(t)\leq1$.

For the quasi-static separation with infinite time $\tau$, the minimum
work is reached as $W_{\mathrm{min}}^{(0)}=-N_{t}k_{B}T_{0}(\epsilon_{\alpha}\ln\epsilon_{\alpha}+\epsilon_{\beta}\ln\epsilon_{\beta})$
by choosing the final volume proportional to the ratio $\epsilon_{\sigma}\equiv N_{\sigma}/N_{t}$
of the two gases, namely, $x_{L}(\tau)=\epsilon_{\alpha}$ and $x_{R}(\tau)=\epsilon_{\beta}$.
See Supplementary Material for detailed discussion. For slow processes
with long operation time $\tau$, the leading term of the excess work
rate is a quadratic form

\begin{equation}
\begin{split}\dot{W}_{\mathrm{ex}}= & \frac{2N_{t}k_{B}T_{0}}{3\gamma}\left[\frac{\epsilon_{\alpha}\dot{x}_{R}}{1-x_{R}}+\frac{\epsilon_{\beta}\dot{x}_{L}}{1-x_{L}}\right]^{2}\\
 & +\frac{N_{t}k_{B}T_{0}V_{t}}{\mathcal{A}}\Big\{\frac{\epsilon_{\alpha}(1-x_{R})}{\mu_{\alpha}}\left[\frac{d}{dt}(\frac{x_{L}}{1-x_{R}})\right]^{2}\\
 & +\frac{\epsilon_{\beta}(1-x_{L})}{\mu_{\beta}}\left[\frac{d}{dt}(\frac{x_{R}}{1-x_{L}})\right]^{2}\Big\}\,.
\end{split}
\label{eq:quadextrawork}
\end{equation}
The first term shows the contribution due to the temperature difference
between the system and the bath during the finite-time separation
process. The second term shows the contribution due to the particle-density
difference across the membranes A and B.

\textit{Riemann geometry of the configuration space.} The quadratic
form of excess work $\dot{W}_{\mathrm{ex}}$ in Eq. (\ref{eq:quadextrawork})
allows the definition of a geometric length in the configuration space
spanned by $(x_{L},x_{R})$. We introduce the timescale of heat transfer
$\tau_{h}\equiv1/\gamma$ and that of particle transport of type-$\sigma$
molecules $\tau_{\sigma}\equiv V_{t}/(\mu_{\sigma}\mathcal{A})$,
and rewrite Eq. (\ref{eq:quadextrawork}) into a compact form 
\begin{align}
\dot{W}_{\mathrm{ex}} & =\left(\begin{array}{cc}
\dot{x}_{L} & \dot{x}_{R}\end{array}\right)G\left(\begin{array}{c}
\dot{x}_{L}\\
\dot{x}_{R}
\end{array}\right)\,,\label{eq:wexdot}
\end{align}
where $G$ is the metric of the current Riemann manifold

\begin{equation}
G=\left(\begin{array}{cc}
g_{LL} & g_{LR}\\
g_{RL} & g_{RR}
\end{array}\right)\,,\label{eq:metric}
\end{equation}
with the components $g_{LL}=N_{t}k_{B}T_{0}\big[\epsilon_{\beta}\tau_{\beta}x_{R}^{2}/(1-x_{L})^{3}+\epsilon_{\alpha}\tau_{\alpha}/(1-x_{R})+2/3\cdot\epsilon_{\beta}^{2}\tau_{h}/(1-x_{L})^{2}\big]\,$,$\quad g_{LR}=g_{RL}=N_{t}k_{B}T_{0}\big[\epsilon_{\alpha}\tau_{\alpha}x_{L}/(1-x_{R})^{2}+\epsilon_{\beta}\tau_{\beta}x_{R}/(1-x_{L})^{2}+2/3\cdot\epsilon_{\alpha}\epsilon_{\beta}\tau_{h}/(1-x_{L})(1-x_{R})\big]\,$,$\quad g_{RR}=N_{t}k_{B}T_{0}\big[\epsilon_{\alpha}\tau_{\alpha}x_{L}^{2}/(1-x_{R})^{3}+\epsilon_{\beta}\tau_{\beta}/(1-x_{L})+2/3\cdot\epsilon_{\alpha}^{2}\tau_{h}/(1-x_{R})^{2}\big]\,$.

The protocol of the separation process is given by $\tilde{x}_{L}(s)=x_{L}(t)$
and $\tilde{x}_{R}(s)=x_{R}(t)$ with the rescaled time $s\equiv t/\tau,\,0\leq s\leq1$.
The protocol $\tilde{x}_{L}(s)$ and $\tilde{x}_{R}(s)$ can be designed
to minimize the excess work with fixed operation time $\tau$ along
a given path $\mathcal{P}(x_{L},x_{R})=0$ in the configuration space.
With the Cauchy-Schwarz inequality $\int_{0}^{\tau}\dot{W}_{\mathrm{ex}}dt\int_{0}^{\tau}dt\geq(\int_{0}^{\tau}\sqrt{\dot{W}_{\mathrm{ex}}}dt)^{2}\equiv\mathcal{L}^{2}$,
the excess work of the given path with different protocols is bounded
by

\begin{equation}
W_{\mathrm{ex}}\ge\frac{\mathcal{L}^{2}}{\tau}\,,\label{eq:bounded}
\end{equation}
where the thermodynamic length $\mathcal{L}=\int_{0}^{\tau}\sqrt{\dot{W}_{\mathrm{ex}}}dt$
is only determined by the path $\mathcal{P}$ in the configuration
space \citep{Salamon1983,Sivak2012}

\begin{equation}
\mathcal{L}=\int_{\mathcal{P}}\sqrt{\left(\begin{array}{cc}
dx_{L} & dx_{R}\end{array}\right)G\left(\begin{array}{c}
dx_{L}\\
dx_{R}
\end{array}\right)}\,.\label{eq:lengthdefine}
\end{equation}
The equality is reached for the protocol with a constant work rate
$\dot{W}_{\mathrm{ex}}=\mathrm{const}$ or a constant velocity of
the thermodynamic length $d\mathcal{L}/ds=\mathrm{const}$. The task
to seek the minimum energy consumption in finite time is converted
into searching geodesic paths connecting the start point $(0,0)$
and the target point $(\epsilon_{\alpha},\epsilon_{\beta})$ in the
configuration space.

Geodesic paths in the current Riemann manifold are described by the
geodesic equations,

\begin{equation}
\frac{d^{2}x^{k}}{dr^{2}}+\Gamma_{\;\,ij}^{k}\frac{dx^{i}}{dr}\frac{dx^{j}}{dr}=0\,,i=1,2\,,\label{eq:geodesic equation}
\end{equation}
where the Einstein notation is employed, and $x^{1}\equiv x_{L},\,x^{2}\equiv x_{R}$.
We parameterize the geodesic paths with the arc length $r$ by setting
the initial condition $g_{ij}\frac{dx^{i}}{dr}\frac{dx^{j}}{dr}=1$.
The Christoffel symbols $\Gamma_{\;\,ij}^{k}$ are obtained as $\Gamma_{\;\,ij}^{k}=\frac{1}{2}g^{kl}\left(\partial g_{li}/\partial x^{j}+\partial g_{lj}/\partial x^{i}-\partial g_{ij}/\partial x^{l}\right),$
where $g^{kl}$ are the elements of the inverse metric $G^{-1}\,$.
The explicit form of $\Gamma_{\;\,ij}^{k}$ is shown in Supplementary
Material.

\textit{Symmetric control path.} For symmetric parameters $\tau_{\alpha}=\tau_{\beta}=\tau_{p},\epsilon_{\alpha}=\epsilon_{\beta}=0.5$,
the symmetric path $x_{L}(t)=x_{R}(t)\equiv x(t)=\tilde{x}(t/\tau)$
is a geodesic path. \begin{widetext}According to Eq. (\ref{eq:lengthdefine}),
the thermodynamic length of the symmetric path as a function of the
endpoint coordinate $x(\tau)$ is obtained as

\begin{equation}
\begin{split}\mathcal{L}_{s}(x(\tau)) & =\int_{0}^{x(\tau)}\sqrt{g_{LL}(x,x)dx^{2}+g_{LR}(x,x)dx^{2}+g_{RL}(x,x)dx^{2}+g_{RR}(x,x)dx^{2}}\,\\
 & =2\sqrt{N_{t}k_{B}T_{0}}\left.\left(\sqrt{\frac{\tau_{p}}{1-x}+\frac{2}{3}\tau_{h}}-\sqrt{\frac{2\tau_{h}}{3}}\sinh^{-1}\sqrt{\frac{2(1-x)\tau_{h}}{3\tau_{p}}}\right)\right|_{x=0}^{x(\tau)}
\end{split}
\label{eq:straightlinelength}
\end{equation}
By setting $x(\tau)=0.5$, we obtain the length of the whole path
$\mathcal{L}_{\mathrm{sym}}=\mathcal{L}_{s}(0.5)\,$.\end{widetext}

The lower bound of excess work for symmetric protocols is $\mathcal{L}_{\mathrm{sym}}^{2}/\tau$.
To reach such bound, we obtain a protocol $\tilde{x}(s)$ in an implicitly
form as
\begin{equation}
\mathcal{L}_{s}(\tilde{x}(s))=s\mathcal{L}_{\mathrm{sym}}.
\end{equation}
It can be verified that $x^{1}(r)=x^{2}(r)=\tilde{x}(r/\mathcal{L}_{\mathrm{sym}})$
is a solution to the geodesic equations.\textbf{ }We explicitly show
two situations where relaxation of the particle transport or that
of heat exchange dominates during the separation process.

(1) Particle transport dominated process ($\tau_{h}\ll\tau_{p}$).
The temperature of the system is identical to that of the bath, $T=T_{0}$.
The excess work rate is simplified as $\dot{W}_{\mathrm{ex}}=N_{t}k_{B}T_{0}\tau_{p}(1-x)\{d/dt[x/(1-x)]\}^{2}$.
The thermodynamic length follows as $\mathcal{L}=(2\sqrt{2}-2)\sqrt{N_{t}k_{B}T_{0}\tau_{p}}$.
According to Eq. (\ref{eq:bounded}), the minimum excess work is

\[
W_{\mathrm{ex}}^{(\mathrm{min})}=(12-8\sqrt{2})N_{t}k_{B}T_{0}\frac{\tau_{p}}{\tau}\,.
\]
The optimal protocol is designed as $\tilde{x}_{L}(s)=\tilde{x}_{R}(s)=1-[(\sqrt{2}-1)s+1]^{-2}$
to achieve the above minimum excess work.

(2) Heat exchange dominated process ($\tau_{h}\gg\tau_{p}$). The
excess work rate is simplified into $\dot{W}_{\mathrm{ex}}=2/3\cdot N_{t}k_{B}T_{0}\tau_{h}[\dot{x}/(1-x)]^{2}$.
The minimum excess work is

\begin{equation}
W_{\mathrm{ex}}^{(\mathrm{min})}=\frac{2(\ln2)^{2}}{3}N_{t}k_{B}T_{0}\frac{\tau_{h}}{\tau}\,.
\end{equation}
The designed protocol is obtained as $\tilde{x}_{L}(s)=\tilde{x}_{R}(s)=1-2^{-s}$.
We summarize the two situations in Table I.

\begin{table*}
\caption{Optimal separation protocol of equally mixed gases under the symmetric
path to move two membranes. Two regions of the relaxation timescales
are considered, (1) particle transport dominated process ($\tau_{h}\ll\tau_{p}$)
and (2) heat exchange dominated process ($\tau_{h}\gg\tau_{p}$).}

\centering{}%
\begin{tabular}{>{\centering}p{3cm}>{\centering}p{8.5cm}>{\centering}p{5cm}}
\toprule 
 & $\tau_{h}\ll\tau_{p}$ & $\tau_{h}\gg\tau_{p}$\tabularnewline
\midrule
\midrule 
$W_{\mathrm{ex}}^{(\mathrm{min})}$ & $(12-8\sqrt{2})N_{t}k_{B}T_{0}\frac{\tau_{p}}{\tau}$ & $\frac{2(\ln2)^{2}}{3}N_{t}k_{B}T_{0}\frac{\tau_{h}}{\tau}.$\tabularnewline
\midrule 
Protocol & $\tilde{x}_{L}(s)=\tilde{x}_{R}(s)=1-[(\sqrt{2}-1)s+1]^{-2}$ & $\tilde{x}_{L}(s)=\tilde{x}_{R}(s)=1-2^{-s}$\tabularnewline
\bottomrule
\end{tabular}
\end{table*}

\begin{figure}[tbph]
\includegraphics{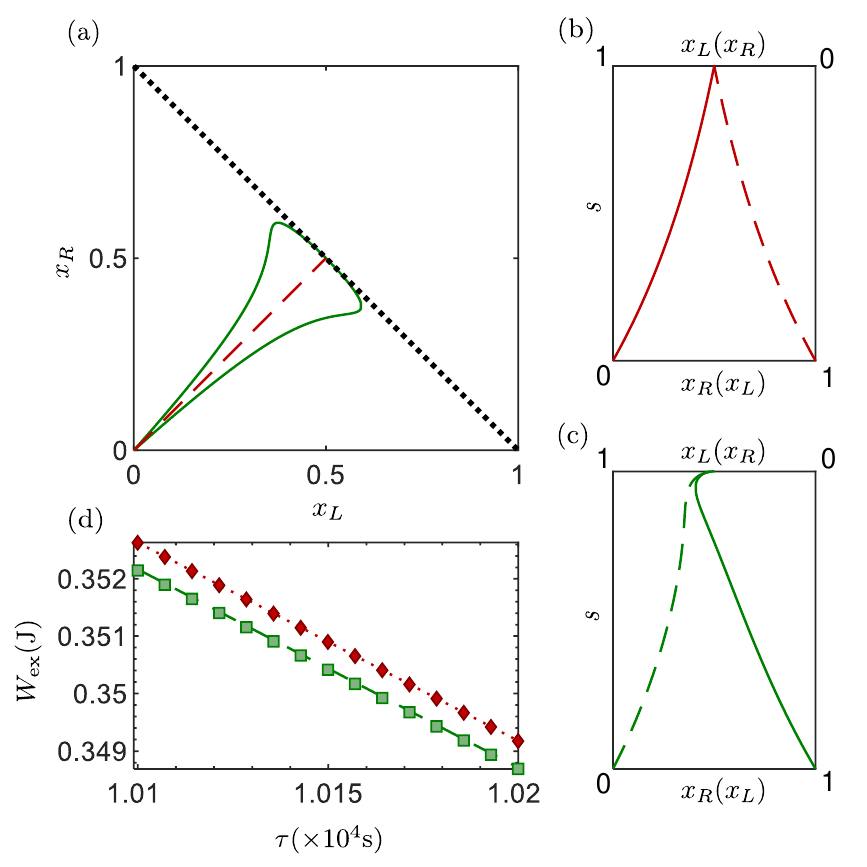}

\caption{Symmetry-breaking protocol of the separation process with the minimum
excess work. The permeable properties of membranes A and B to the
corresponding type of molecules are the same, and the numbers of type-$\alpha$
and type-$\beta$ molecules are equal. Namely, all the parameters
of this system are symmetric ($\tau_{\alpha}=\tau_{\beta}=1\mathrm{s},\,\tau_{h}=0.1\mathrm{s},\:\epsilon_{\alpha}=\epsilon_{\beta}=0.5,N_{t}=2N_{A},T_{0}=298.15\mathrm{K}$,
where $N_{A}$ is the Avogadro constant). (a) Three geodesic paths
connecting $(0,0)$ and $(0.5,0.5)$. The thermodynamic length of
the two symmetry-breaking geodesic paths (green solid line) is $\mathcal{L}_{\mathrm{geo}}=59.64(\mathrm{J\cdot s})^{1/2}$,
which is shorter than that of the symmetric one (red dashed line)
$\mathcal{L}_{\mathrm{sym}}=59.68(\mathrm{J\cdot s})^{1/2}$. (b)
and (c) Control schemes of the protocols of the symmetry-breaking
and the symmetric paths. In the symmetry-breaking protocol, the two
membranes approach each other at a position slightly different from
the equilibrium position, and are moved together to the equilibrium
position. (d) Excess work $W_{\mathrm{ex}}$ versus operation time
$\tau\,$. The red diamonds (green squares) represent the irreversible
work obtained by solving the evolution equation of the system along
the symmetric (symmetry-breaking) path. While the red dashed line
(green solid line) represents the excess work predicted by the thermodynamic
length $W_{\mathrm{ex}}=\mathcal{L}^{2}/\tau$.}

\label{fig:The-optimal-protocol}
\end{figure}

\textit{Symmetry breaking in the optimal separation protocol.} The
question arises that whether the straightforward symmetric protocol
above is the optimal one with the minimum energy expense. Our answer
is no. We find a symmetry-breaking protocol for the perfect symmetric
setup to reach the minimum energy expense.

For given parameters, we use the shooting method \citep{Berger2007}
to numerically solve the geodesic path connecting the initial position
$(0,0)$ and the final position $(0.5,0.5)$ in the configuration
space. For room temperature $T_{0}=298.15\mathrm{K}$, we consider
the symmetric situation of two moles equally mixed gases with $N_{t}=2N_{A}$,
where $N_{A}$ is the Avogadro constant. The particle transport relaxation
timescales are the same $\tau_{\alpha}=\tau_{\beta}=1\mathrm{s}$,
and the heat conduction relaxation timescale is $\tau_{h}=0.1\mathrm{s}$.
We find three geodesic paths connecting the start point $(0,0)$ and
the target point $(0.5,0.5)$, shown in Fig \ref{fig:The-optimal-protocol}(a)
as the symmetric red dashed line and the two green lines. To visualize
the symmetry breaking of the control scheme in the optimal separation
protocol, we sketch the three-dimensional embedding of the current
two-dimensional Riemann manifold in Fig. \ref{fig:scheme}(b), where
the thermodynamic length is reflected by the length of the path. The
symmetric and the symmetry-breaking geodesic paths are shown in red
and green lines. The green lines are shorter than the red one, representing
smaller excess work.

The red dashed line of symmetric protocol shown in Fig \ref{fig:The-optimal-protocol}(a)
has a thermodynamic length of $\mathcal{L}_{\mathrm{sym}}=59.68(\mathrm{J\cdot s})^{1/2}$.
The lengths of the symmetry-breaking green geodesic paths are $\mathcal{L}_{\mathrm{geo}}=59.64(\mathrm{J\cdot s})^{1/2}$,
shorter than the symmetric one. The control schemes $\tilde{x}_{L,R}(s)$
for the symmetric and symmetry-breaking protocol are shown in Fig.
\ref{fig:The-optimal-protocol}(b) and Fig. \ref{fig:The-optimal-protocol}(c).
In the symmetry-breaking protocol shown in Fig. \ref{fig:The-optimal-protocol}(c),
the two membranes approach each other at a position different from
the equilibrium position, and then are moved together to the equilibrium
position. In Fig. \ref{fig:The-optimal-protocol}(d), we compares
the results obtained with the quadratic form in Eq. (\ref{eq:quadextrawork})
and the numerical solution of the evolution equations (\ref{eq:particle transport})
and (\ref{eq:temperature}). The red dotted line and the green dashed
line show the approximate excess work predicted by thermodynamic length
$W_{\mathrm{ex}}=\mathcal{L}^{2}/\tau$. And the red diamonds and
the green squares present the excess work calculated with Eq. (\ref{eq:W_accurate})
by solving the evolution Eqs. (\ref{eq:particle transport}) and (\ref{eq:temperature})
for the corresponding protocols $V_{L,R}(t)$ obtained from the geodesic
equations.

\textit{Conclusion.} We prove the equivalence between designing optimum
control to achieve the minimum energy expense and finding the geodesic
path in a geometric space. Such an equivalence has also been exploited
in the optimization of control protocols in stochastic and quantum
thermodynamics \citep{Sivak2012,Scandi2018,Chen2021,Li2022}. With
this equivalence, we show the minimum excess energy expense is proportional
to the square of the length of that geodesic path and inversely proportional
to the operation time $\tau$. In a separation process where all the
parameters are symmetric $\tau_{\alpha}=\tau_{\beta}=\tau_{p}\,,\epsilon_{\alpha}=\epsilon_{\beta}=0.5\,$,
we found three geodesic paths, among them a simple and straightforward
path is symmetric in the configuration space $(x_{L},x_{R})$. The
corresponding protocol is to move the two membranes with the same
speed. In this situation, we predict that the complete separation
of equally mixed single-atom gases requires excess work at least $(12-8\sqrt{2})N_{t}k_{B}T_{0}\tau_{p}/\tau$
for the particle transport dominated process and $2(\ln2)^{2}N_{t}k_{B}T_{0}\tau_{h}/(3\tau)$
for the heat exchange dominated process, where $\tau_{p}$ and $\tau_{h}$
are the relaxation time for the particle transport and the heat exchange.

Importantly, the current work shows that such symmetric protocol is
not the best one to achieve the minimal energy expense. A symmetry-breaking
protocol is the optimal protocol that achieves the lowest energy expense
for given operation time $\tau$. In the optimal protocol, one membrane
is moved faster than the other, the membranes approach each other
at a position slightly deviated from the equilibrium position. Then
they are moved together to the equilibrium position.

\textit{Acknowledgments.} This work is supported by the National Natural
Science Foundation of China (NSFC) (Grants No. 12088101, No. 11534002,
No. 11875049, No. U1930402, No. U1930403 and No. 12047549) and the
National Basic Research Program of China (Grant No. 2016YFA0301201).

Jin-Fu Chen and Ruo-Xun Zhai contributed equally to this work.

\bibliography{separationbib}

\begin{thebibliography}{34}%
\makeatletter
\providecommand \@ifxundefined [1]{%
 \@ifx{#1\undefined}
}%
\providecommand \@ifnum [1]{%
 \ifnum #1\expandafter \@firstoftwo
 \else \expandafter \@secondoftwo
 \fi
}%
\providecommand \@ifx [1]{%
 \ifx #1\expandafter \@firstoftwo
 \else \expandafter \@secondoftwo
 \fi
}%
\providecommand \natexlab [1]{#1}%
\providecommand \enquote  [1]{``#1''}%
\providecommand \bibnamefont  [1]{#1}%
\providecommand \bibfnamefont [1]{#1}%
\providecommand \citenamefont [1]{#1}%
\providecommand \href@noop [0]{\@secondoftwo}%
\providecommand \href [0]{\begingroup \@sanitize@url \@href}%
\providecommand \@href[1]{\@@startlink{#1}\@@href}%
\providecommand \@@href[1]{\endgroup#1\@@endlink}%
\providecommand \@sanitize@url [0]{\catcode `\\12\catcode `\$12\catcode
  `\&12\catcode `\#12\catcode `\^12\catcode `\_12\catcode `\%12\relax}%
\providecommand \@@startlink[1]{}%
\providecommand \@@endlink[0]{}%
\providecommand \url  [0]{\begingroup\@sanitize@url \@url }%
\providecommand \@url [1]{\endgroup\@href {#1}{\urlprefix }}%
\providecommand \urlprefix  [0]{URL }%
\providecommand \Eprint [0]{\href }%
\providecommand \doibase [0]{http://dx.doi.org/}%
\providecommand \selectlanguage [0]{\@gobble}%
\providecommand \bibinfo  [0]{\@secondoftwo}%
\providecommand \bibfield  [0]{\@secondoftwo}%
\providecommand \translation [1]{[#1]}%
\providecommand \BibitemOpen [0]{}%
\providecommand \bibitemStop [0]{}%
\providecommand \bibitemNoStop [0]{.\EOS\space}%
\providecommand \EOS [0]{\spacefactor3000\relax}%
\providecommand \BibitemShut  [1]{\csname bibitem#1\endcsname}%
\let\auto@bib@innerbib\@empty
\bibitem [{\citenamefont {Sholl}\ and\ \citenamefont
  {Lively}(2016)}]{Sholl2016}%
  \BibitemOpen
  \bibfield  {author} {\bibinfo {author} {\bibfnamefont {David~S.}\
  \bibnamefont {Sholl}}\ and\ \bibinfo {author} {\bibfnamefont {Ryan~P.}\
  \bibnamefont {Lively}},\ }\bibfield  {title} {\enquote {\bibinfo {title}
  {Seven chemical separations to change the world},}\ }\href {\doibase
  10.1038/532435a} {\bibfield  {journal} {\bibinfo  {journal} {Nature}\
  }\textbf {\bibinfo {volume} {532}},\ \bibinfo {pages} {435--437} (\bibinfo
  {year} {2016})}\BibitemShut {NoStop}%
\bibitem [{\citenamefont {Koros}\ and\ \citenamefont
  {Lively}(2012)}]{Koros2012}%
  \BibitemOpen
  \bibfield  {author} {\bibinfo {author} {\bibfnamefont {William~J.}\
  \bibnamefont {Koros}}\ and\ \bibinfo {author} {\bibfnamefont {Ryan~P.}\
  \bibnamefont {Lively}},\ }\bibfield  {title} {\enquote {\bibinfo {title}
  {Water and beyond: Expanding the spectrum of large-scale energy efficient
  separation processes},}\ }\href {\doibase https://doi.org/10.1002/aic.13888}
  {\bibfield  {journal} {\bibinfo  {journal} {AIChE Journal}\ }\textbf
  {\bibinfo {volume} {58}},\ \bibinfo {pages} {2624--2633} (\bibinfo {year}
  {2012})}\BibitemShut {NoStop}%
\bibitem [{\citenamefont {Alvarez}\ \emph {et~al.}(2018)\citenamefont
  {Alvarez}, \citenamefont {Chan}, \citenamefont {Elimelech}, \citenamefont
  {Halas},\ and\ \citenamefont {Villagr{\'{a}}n}}]{Alvarez2018}%
  \BibitemOpen
  \bibfield  {author} {\bibinfo {author} {\bibfnamefont {Pedro J.~J.}\
  \bibnamefont {Alvarez}}, \bibinfo {author} {\bibfnamefont {Candace~K.}\
  \bibnamefont {Chan}}, \bibinfo {author} {\bibfnamefont {Menachem}\
  \bibnamefont {Elimelech}}, \bibinfo {author} {\bibfnamefont {Naomi~J.}\
  \bibnamefont {Halas}}, \ and\ \bibinfo {author} {\bibfnamefont {Dino}\
  \bibnamefont {Villagr{\'{a}}n}},\ }\bibfield  {title} {\enquote {\bibinfo
  {title} {Emerging opportunities for nanotechnology to enhance water
  security},}\ }\href {\doibase 10.1038/s41565-018-0203-2} {\bibfield
  {journal} {\bibinfo  {journal} {Nat. Nanotechnol.}\ }\textbf {\bibinfo
  {volume} {13}},\ \bibinfo {pages} {634--641} (\bibinfo {year}
  {2018})}\BibitemShut {NoStop}%
\bibitem [{\citenamefont {Noamani}\ \emph {et~al.}(2019)\citenamefont
  {Noamani}, \citenamefont {Niroomand}, \citenamefont {Rastgar},\ and\
  \citenamefont {Sadrzadeh}}]{Noamani2019}%
  \BibitemOpen
  \bibfield  {author} {\bibinfo {author} {\bibfnamefont {Sadaf}\ \bibnamefont
  {Noamani}}, \bibinfo {author} {\bibfnamefont {Shirin}\ \bibnamefont
  {Niroomand}}, \bibinfo {author} {\bibfnamefont {Masoud}\ \bibnamefont
  {Rastgar}}, \ and\ \bibinfo {author} {\bibfnamefont {Mohtada}\ \bibnamefont
  {Sadrzadeh}},\ }\bibfield  {title} {\enquote {\bibinfo {title} {Carbon-based
  polymer nanocomposite membranes for oily wastewater treatment},}\ }\href
  {\doibase 10.1038/s41545-019-0044-z} {\bibfield  {journal} {\bibinfo
  {journal} {npj Clean Water}\ }\textbf {\bibinfo {volume} {2}},\ \bibinfo
  {pages} {20} (\bibinfo {year} {2019})}\BibitemShut {NoStop}%
\bibitem [{\citenamefont {van Reis}\ and\ \citenamefont
  {Zydney}(2007)}]{Reis2007}%
  \BibitemOpen
  \bibfield  {author} {\bibinfo {author} {\bibfnamefont {Robert}\ \bibnamefont
  {van Reis}}\ and\ \bibinfo {author} {\bibfnamefont {Andrew}\ \bibnamefont
  {Zydney}},\ }\bibfield  {title} {\enquote {\bibinfo {title} {Bioprocess
  membrane technology},}\ }\href {\doibase 10.1016/j.memsci.2007.02.045}
  {\bibfield  {journal} {\bibinfo  {journal} {J. Membr. Sci.}\ }\textbf
  {\bibinfo {volume} {297}},\ \bibinfo {pages} {16--50} (\bibinfo {year}
  {2007})}\BibitemShut {NoStop}%
\bibitem [{\citenamefont {Xie}\ \emph {et~al.}(2008)\citenamefont {Xie},
  \citenamefont {Chu},\ and\ \citenamefont {Deng}}]{Xie2008}%
  \BibitemOpen
  \bibfield  {author} {\bibinfo {author} {\bibfnamefont {Rui}\ \bibnamefont
  {Xie}}, \bibinfo {author} {\bibfnamefont {Liang-Yin}\ \bibnamefont {Chu}}, \
  and\ \bibinfo {author} {\bibfnamefont {Jin-Gen}\ \bibnamefont {Deng}},\
  }\bibfield  {title} {\enquote {\bibinfo {title} {Membranes and membrane
  processes for chiral resolution},}\ }\href {\doibase 10.1039/b713350b}
  {\bibfield  {journal} {\bibinfo  {journal} {Chem. Soc. Rev.}\ }\textbf
  {\bibinfo {volume} {37}},\ \bibinfo {pages} {1243} (\bibinfo {year}
  {2008})}\BibitemShut {NoStop}%
\bibitem [{\citenamefont {O.~Falk-Pedersen}(1997)}]{FalkPedersen1997}%
  \BibitemOpen
  \bibfield  {author} {\bibinfo {author} {\bibfnamefont {H.~Dannstr\"om}\
  \bibnamefont {O.~Falk-Pedersen}},\ }\bibfield  {title} {\enquote {\bibinfo
  {title} {Separation of carbon dioxide from offshore gas turbine exhaust},}\
  }\href {\doibase 10.1016/s0196-8904(96)00250-6} {\bibfield  {journal}
  {\bibinfo  {journal} {Energy Conversion and Management}\ }\textbf {\bibinfo
  {volume} {38}},\ \bibinfo {pages} {S81--S86} (\bibinfo {year}
  {1997})}\BibitemShut {NoStop}%
\bibitem [{\citenamefont {Favre}(2007)}]{Favre2007}%
  \BibitemOpen
  \bibfield  {author} {\bibinfo {author} {\bibfnamefont {Eric}\ \bibnamefont
  {Favre}},\ }\bibfield  {title} {\enquote {\bibinfo {title} {Carbon dioxide
  recovery from post-combustion processes: Can gas permeation membranes compete
  with absorption?}}\ }\href {\doibase 10.1016/j.memsci.2007.02.007} {\bibfield
   {journal} {\bibinfo  {journal} {Journal of Membrane Science}\ }\textbf
  {\bibinfo {volume} {294}},\ \bibinfo {pages} {50--59} (\bibinfo {year}
  {2007})}\BibitemShut {NoStop}%
\bibitem [{\citenamefont {Merkel}\ \emph {et~al.}(2010)\citenamefont {Merkel},
  \citenamefont {Lin}, \citenamefont {Wei},\ and\ \citenamefont
  {Baker}}]{Merkel2010}%
  \BibitemOpen
  \bibfield  {author} {\bibinfo {author} {\bibfnamefont {Tim~C.}\ \bibnamefont
  {Merkel}}, \bibinfo {author} {\bibfnamefont {Haiqing}\ \bibnamefont {Lin}},
  \bibinfo {author} {\bibfnamefont {Xiaotong}\ \bibnamefont {Wei}}, \ and\
  \bibinfo {author} {\bibfnamefont {Richard}\ \bibnamefont {Baker}},\
  }\bibfield  {title} {\enquote {\bibinfo {title} {Power plant post-combustion
  carbon dioxide capture: An opportunity for membranes},}\ }\href {\doibase
  10.1016/j.memsci.2009.10.041} {\bibfield  {journal} {\bibinfo  {journal}
  {Journal of Membrane Science}\ }\textbf {\bibinfo {volume} {359}},\ \bibinfo
  {pages} {126--139} (\bibinfo {year} {2010})}\BibitemShut {NoStop}%
\bibitem [{\citenamefont {Adewole}\ \emph {et~al.}(2013)\citenamefont
  {Adewole}, \citenamefont {Ahmad}, \citenamefont {Ismail},\ and\ \citenamefont
  {Leo}}]{Adewole2013}%
  \BibitemOpen
  \bibfield  {author} {\bibinfo {author} {\bibfnamefont {J.K.}\ \bibnamefont
  {Adewole}}, \bibinfo {author} {\bibfnamefont {A.L.}\ \bibnamefont {Ahmad}},
  \bibinfo {author} {\bibfnamefont {S.}~\bibnamefont {Ismail}}, \ and\ \bibinfo
  {author} {\bibfnamefont {C.P.}\ \bibnamefont {Leo}},\ }\bibfield  {title}
  {\enquote {\bibinfo {title} {Current challenges in membrane separation of
  {CO}2 from natural gas: A review},}\ }\href {\doibase
  10.1016/j.ijggc.2013.04.012} {\bibfield  {journal} {\bibinfo  {journal}
  {International Journal of Greenhouse Gas Control}\ }\textbf {\bibinfo
  {volume} {17}},\ \bibinfo {pages} {46--65} (\bibinfo {year}
  {2013})}\BibitemShut {NoStop}%
\bibitem [{\citenamefont {Britan}\ \emph {et~al.}(1991)\citenamefont {Britan},
  \citenamefont {Leites},\ and\ \citenamefont {Vasilkovskaya}}]{Britan1991}%
  \BibitemOpen
  \bibfield  {author} {\bibinfo {author} {\bibfnamefont {I.M.}\ \bibnamefont
  {Britan}}, \bibinfo {author} {\bibfnamefont {I.L.}\ \bibnamefont {Leites}}, \
  and\ \bibinfo {author} {\bibfnamefont {T.N.}\ \bibnamefont {Vasilkovskaya}},\
  }\bibfield  {title} {\enquote {\bibinfo {title} {Membrane technology of
  mixed-gas separation: thermodynamic analysis for feasibility study},}\ }\href
  {\doibase 10.1016/s0376-7388(00)80588-8} {\bibfield  {journal} {\bibinfo
  {journal} {J. Membr. Sci.}\ }\textbf {\bibinfo {volume} {55}},\ \bibinfo
  {pages} {349--352} (\bibinfo {year} {1991})}\BibitemShut {NoStop}%
\bibitem [{\citenamefont {Koros}\ and\ \citenamefont
  {Fleming}(1993)}]{Koros1993}%
  \BibitemOpen
  \bibfield  {author} {\bibinfo {author} {\bibfnamefont {W.J.}\ \bibnamefont
  {Koros}}\ and\ \bibinfo {author} {\bibfnamefont {G.K.}\ \bibnamefont
  {Fleming}},\ }\bibfield  {title} {\enquote {\bibinfo {title} {Membrane-based
  gas separation},}\ }\href {\doibase 10.1016/0376-7388(93)80013-n} {\bibfield
  {journal} {\bibinfo  {journal} {J. Membr. Sci.}\ }\textbf {\bibinfo {volume}
  {83}},\ \bibinfo {pages} {1--80} (\bibinfo {year} {1993})}\BibitemShut
  {NoStop}%
\bibitem [{\citenamefont {Spillman}(1995)}]{Spillman1995}%
  \BibitemOpen
  \bibfield  {author} {\bibinfo {author} {\bibfnamefont {Robert}\ \bibnamefont
  {Spillman}},\ }\bibfield  {title} {\enquote {\bibinfo {title} {Chapter 13
  economics of gas separation membrane processes},}\ }in\ \href {\doibase
  10.1016/s0927-5193(06)80015-x} {\emph {\bibinfo {booktitle} {Membrane Science
  and Technology}}}\ (\bibinfo  {publisher} {Elsevier},\ \bibinfo {year}
  {1995})\ pp.\ \bibinfo {pages} {589--667}\BibitemShut {NoStop}%
\bibitem [{\citenamefont {Pabby}\ \emph {et~al.}(2015)\citenamefont {Pabby},
  \citenamefont {Rizvi},\ and\ \citenamefont {Requena}}]{Pabby2015}%
  \BibitemOpen
  \bibfield  {author} {\bibinfo {author} {\bibfnamefont {Anil~Kumar}\
  \bibnamefont {Pabby}}, \bibinfo {author} {\bibfnamefont {Syed S.~H.}\
  \bibnamefont {Rizvi}}, \ and\ \bibinfo {author} {\bibfnamefont {Ana
  Maria~Sastre}\ \bibnamefont {Requena}},\ }\href@noop {} {\emph {\bibinfo
  {title} {Handbook of Membrane Separations Chemical, Pharmaceutical, Food, and
  Biotechnological Applications, Second Edition}}}\ (\bibinfo  {publisher}
  {Taylor and Francis Group},\ \bibinfo {year} {2015})\ p.\ \bibinfo {pages}
  {878}\BibitemShut {NoStop}%
\bibitem [{\citenamefont {Castel}\ and\ \citenamefont
  {Favre}(2018)}]{Castel2018}%
  \BibitemOpen
  \bibfield  {author} {\bibinfo {author} {\bibfnamefont {Christophe}\
  \bibnamefont {Castel}}\ and\ \bibinfo {author} {\bibfnamefont {Eric}\
  \bibnamefont {Favre}},\ }\bibfield  {title} {\enquote {\bibinfo {title}
  {Membrane separations and energy efficiency},}\ }\href {\doibase
  10.1016/j.memsci.2017.11.035} {\bibfield  {journal} {\bibinfo  {journal} {J.
  Membr. Sci.}\ }\textbf {\bibinfo {volume} {548}},\ \bibinfo {pages}
  {345--357} (\bibinfo {year} {2018})}\BibitemShut {NoStop}%
\bibitem [{\citenamefont {Purkait}\ and\ \citenamefont
  {Singh}(2018)}]{Purkait2018}%
  \BibitemOpen
  \bibfield  {author} {\bibinfo {author} {\bibfnamefont {Mihir~K.}\
  \bibnamefont {Purkait}}\ and\ \bibinfo {author} {\bibfnamefont {Randeep}\
  \bibnamefont {Singh}},\ }\href@noop {} {\emph {\bibinfo {title} {Membrane
  Technology in Separation Science}}}\ (\bibinfo  {publisher} {Taylor and
  Francis Group},\ \bibinfo {year} {2018})\ p.\ \bibinfo {pages}
  {242}\BibitemShut {NoStop}%
\bibitem [{\citenamefont {Weinhold}(1975)}]{Weinhold1975}%
  \BibitemOpen
  \bibfield  {author} {\bibinfo {author} {\bibfnamefont {F.}~\bibnamefont
  {Weinhold}},\ }\bibfield  {title} {\enquote {\bibinfo {title} {Metric
  geometry of equilibrium thermodynamics},}\ }\href {\doibase 10.1063/1.431689}
  {\bibfield  {journal} {\bibinfo  {journal} {J. Chem. Phys.}\ }\textbf
  {\bibinfo {volume} {63}},\ \bibinfo {pages} {2479--2483} (\bibinfo {year}
  {1975})}\BibitemShut {NoStop}%
\bibitem [{\citenamefont {Ruppeiner}(1979)}]{Ruppeiner1979}%
  \BibitemOpen
  \bibfield  {author} {\bibinfo {author} {\bibfnamefont {George}\ \bibnamefont
  {Ruppeiner}},\ }\bibfield  {title} {\enquote {\bibinfo {title}
  {Thermodynamics: A riemannian geometric model},}\ }\href {\doibase
  10.1103/physreva.20.1608} {\bibfield  {journal} {\bibinfo  {journal} {Phys.
  Rev. A}\ }\textbf {\bibinfo {volume} {20}},\ \bibinfo {pages} {1608--1613}
  (\bibinfo {year} {1979})}\BibitemShut {NoStop}%
\bibitem [{\citenamefont {Salamon}\ and\ \citenamefont
  {Berry}(1983)}]{Salamon1983}%
  \BibitemOpen
  \bibfield  {author} {\bibinfo {author} {\bibfnamefont {Peter}\ \bibnamefont
  {Salamon}}\ and\ \bibinfo {author} {\bibfnamefont {R.~Stephen}\ \bibnamefont
  {Berry}},\ }\bibfield  {title} {\enquote {\bibinfo {title} {Thermodynamic
  length and dissipated availability},}\ }\href {\doibase
  10.1103/physrevlett.51.1127} {\bibfield  {journal} {\bibinfo  {journal}
  {Phys. Rev. Lett.}\ }\textbf {\bibinfo {volume} {51}},\ \bibinfo {pages}
  {1127--1130} (\bibinfo {year} {1983})}\BibitemShut {NoStop}%
\bibitem [{\citenamefont {Crooks}(2007)}]{Crooks2007}%
  \BibitemOpen
  \bibfield  {author} {\bibinfo {author} {\bibfnamefont {Gavin~E.}\
  \bibnamefont {Crooks}},\ }\bibfield  {title} {\enquote {\bibinfo {title}
  {Measuring thermodynamic length},}\ }\href {\doibase
  10.1103/physrevlett.99.100602} {\bibfield  {journal} {\bibinfo  {journal}
  {Phys. Rev. Lett.}\ }\textbf {\bibinfo {volume} {99}},\ \bibinfo {pages}
  {100602} (\bibinfo {year} {2007})}\BibitemShut {NoStop}%
\bibitem [{\citenamefont {Sivak}\ and\ \citenamefont
  {Crooks}(2012)}]{Sivak2012}%
  \BibitemOpen
  \bibfield  {author} {\bibinfo {author} {\bibfnamefont {David~A.}\
  \bibnamefont {Sivak}}\ and\ \bibinfo {author} {\bibfnamefont {Gavin~E.}\
  \bibnamefont {Crooks}},\ }\bibfield  {title} {\enquote {\bibinfo {title}
  {Thermodynamic metrics and optimal paths},}\ }\href {\doibase
  10.1103/physrevlett.108.190602} {\bibfield  {journal} {\bibinfo  {journal}
  {Phys. Rev. Lett.}\ }\textbf {\bibinfo {volume} {108}},\ \bibinfo {pages}
  {190602} (\bibinfo {year} {2012})}\BibitemShut {NoStop}%
\bibitem [{\citenamefont {Scandi}\ and\ \citenamefont
  {Perarnau-Llobet}(2019)}]{Scandi2018}%
  \BibitemOpen
  \bibfield  {author} {\bibinfo {author} {\bibfnamefont {Matteo}\ \bibnamefont
  {Scandi}}\ and\ \bibinfo {author} {\bibfnamefont {Mart{\'{\i}}}\ \bibnamefont
  {Perarnau-Llobet}},\ }\bibfield  {title} {\enquote {\bibinfo {title}
  {Thermodynamic length in open quantum systems},}\ }\href {\doibase
  10.22331/q-2019-10-24-197} {\bibfield  {journal} {\bibinfo  {journal}
  {Quantum}\ }\textbf {\bibinfo {volume} {3}},\ \bibinfo {pages} {197}
  (\bibinfo {year} {2019})}\BibitemShut {NoStop}%
\bibitem [{\citenamefont {Chen}\ \emph {et~al.}(2021)\citenamefont {Chen},
  \citenamefont {Sun},\ and\ \citenamefont {Dong}}]{Chen2021}%
  \BibitemOpen
  \bibfield  {author} {\bibinfo {author} {\bibfnamefont {Jin-Fu}\ \bibnamefont
  {Chen}}, \bibinfo {author} {\bibfnamefont {C.~P.}\ \bibnamefont {Sun}}, \
  and\ \bibinfo {author} {\bibfnamefont {Hui}\ \bibnamefont {Dong}},\
  }\bibfield  {title} {\enquote {\bibinfo {title} {Extrapolating the
  thermodynamic length with finite-time measurements},}\ }\href {\doibase
  10.1103/physreve.104.034117} {\bibfield  {journal} {\bibinfo  {journal}
  {Phys. Rev. E}\ }\textbf {\bibinfo {volume} {104}},\ \bibinfo {pages}
  {034117} (\bibinfo {year} {2021})}\BibitemShut {NoStop}%
\bibitem [{\citenamefont {Li}\ \emph {et~al.}(2022)\citenamefont {Li},
  \citenamefont {Chen}, \citenamefont {Sun},\ and\ \citenamefont
  {Dong}}]{Li2022}%
  \BibitemOpen
  \bibfield  {author} {\bibinfo {author} {\bibfnamefont {Geng}\ \bibnamefont
  {Li}}, \bibinfo {author} {\bibfnamefont {Jin-Fu}\ \bibnamefont {Chen}},
  \bibinfo {author} {\bibfnamefont {C.P.}\ \bibnamefont {Sun}}, \ and\ \bibinfo
  {author} {\bibfnamefont {Hui}\ \bibnamefont {Dong}},\ }\bibfield  {title}
  {\enquote {\bibinfo {title} {Geodesic path for the minimal energy cost in
  shortcuts to isothermality},}\ }\href
  {https://doi.org/10.1103/PhysRevLett.128.230603} {\bibfield  {journal}
  {\bibinfo  {journal} {Phys. Rev. Lett.}\ }\textbf {\bibinfo {volume} {128}},\
  \bibinfo {pages} {230603} (\bibinfo {year} {2022})}\BibitemShut {NoStop}%
\bibitem [{\citenamefont {Callen}(1985)}]{book:18204}%
  \BibitemOpen
  \bibfield  {author} {\bibinfo {author} {\bibfnamefont {Herbert~B.}\
  \bibnamefont {Callen}},\ }\href@noop {} {\emph {\bibinfo {title}
  {{T}hermodynamics and an {I}ntroduction to {T}hermostatistics}}},\ \bibinfo
  {edition} {2nd}\ ed.\ (\bibinfo  {publisher} {Wiley},\ \bibinfo {year}
  {1985})\BibitemShut {NoStop}%
\bibitem [{\citenamefont {Huang}(1987)}]{book:18254}%
  \BibitemOpen
  \bibfield  {author} {\bibinfo {author} {\bibfnamefont {Kerson}\ \bibnamefont
  {Huang}},\ }\href
  {http://gen.lib.rus.ec/book/index.php?md5=8BDCF0D4A028500E98805FE0693B6718}
  {\emph {\bibinfo {title} {Statistical mechanics}}},\ \bibinfo {edition}
  {2nd}\ ed.\ (\bibinfo  {publisher} {Wiley},\ \bibinfo {year}
  {1987})\BibitemShut {NoStop}%
\bibitem [{\citenamefont {Tsirlin}\ \emph {et~al.}(2002)\citenamefont
  {Tsirlin}, \citenamefont {Kazakov},\ and\ \citenamefont
  {Zubov}}]{Tsirlin2002}%
  \BibitemOpen
  \bibfield  {author} {\bibinfo {author} {\bibfnamefont {Anatoliy~M.}\
  \bibnamefont {Tsirlin}}, \bibinfo {author} {\bibfnamefont {Vladimir}\
  \bibnamefont {Kazakov}}, \ and\ \bibinfo {author} {\bibfnamefont
  {Dmitrii~V.}\ \bibnamefont {Zubov}},\ }\bibfield  {title} {\enquote {\bibinfo
  {title} {Finite-time thermodynamics:~ limiting possibilities of irreversible
  separation processes{\textdagger}},}\ }\href {\doibase 10.1021/jp025524v}
  {\bibfield  {journal} {\bibinfo  {journal} {J. Phys. Chem. A}\ }\textbf
  {\bibinfo {volume} {106}},\ \bibinfo {pages} {10926--10936} (\bibinfo {year}
  {2002})}\BibitemShut {NoStop}%
\bibitem [{\citenamefont {Xu}\ and\ \citenamefont {Agrawal}(1996)}]{Xu1996}%
  \BibitemOpen
  \bibfield  {author} {\bibinfo {author} {\bibfnamefont {Jianguo}\ \bibnamefont
  {Xu}}\ and\ \bibinfo {author} {\bibfnamefont {Rakesh}\ \bibnamefont
  {Agrawal}},\ }\bibfield  {title} {\enquote {\bibinfo {title} {Membrane
  separation process analysis and design strategies based on thermodynamic
  efficiency of permeation},}\ }\href {\doibase 10.1016/0009-2509(95)00262-6}
  {\bibfield  {journal} {\bibinfo  {journal} {Chem. Eng. Sci.}\ }\textbf
  {\bibinfo {volume} {51}},\ \bibinfo {pages} {365--385} (\bibinfo {year}
  {1996})}\BibitemShut {NoStop}%
\bibitem [{\citenamefont {Sieniutycz}\ and\ \citenamefont
  {Tsirlin}(2017)}]{Sieniutycz2017}%
  \BibitemOpen
  \bibfield  {author} {\bibinfo {author} {\bibfnamefont {Stanislaw}\
  \bibnamefont {Sieniutycz}}\ and\ \bibinfo {author} {\bibfnamefont {Anatoly}\
  \bibnamefont {Tsirlin}},\ }\bibfield  {title} {\enquote {\bibinfo {title}
  {Finding limiting possibilities of thermodynamic systems by optimization},}\
  }\href {\doibase 10.1098/rsta.2016.0219} {\bibfield  {journal} {\bibinfo
  {journal} {Phil. trans. R. Soc. A}\ }\textbf {\bibinfo {volume} {375}},\
  \bibinfo {pages} {20160219} (\bibinfo {year} {2017})}\BibitemShut {NoStop}%
\bibitem [{\citenamefont {Sieniutycz}\ and\ \citenamefont
  {Je{\.{z}}owski}(2018)}]{Sieniutycz2018}%
  \BibitemOpen
  \bibfield  {author} {\bibinfo {author} {\bibfnamefont {Stanis{\l}aw}\
  \bibnamefont {Sieniutycz}}\ and\ \bibinfo {author} {\bibfnamefont {Jacek}\
  \bibnamefont {Je{\.{z}}owski}},\ }\bibfield  {title} {\enquote {\bibinfo
  {title} {Optimization and qualitative aspects of separation systems},}\ }in\
  \href {\doibase 10.1016/b978-0-08-102557-4.00008-6} {\emph {\bibinfo
  {booktitle} {Energy Optimization in Process Systems and Fuel Cells}}}\
  (\bibinfo  {publisher} {Elsevier},\ \bibinfo {year} {2018})\ pp.\ \bibinfo
  {pages} {273--333}\BibitemShut {NoStop}%
\bibitem [{\citenamefont {Curzon}\ and\ \citenamefont
  {Ahlborn}(1975)}]{Curzon1975}%
  \BibitemOpen
  \bibfield  {author} {\bibinfo {author} {\bibfnamefont {F.~L.}\ \bibnamefont
  {Curzon}}\ and\ \bibinfo {author} {\bibfnamefont {B.}~\bibnamefont
  {Ahlborn}},\ }\bibfield  {title} {\enquote {\bibinfo {title} {Efficiency of a
  carnot engine at maximum power output},}\ }\href {\doibase 10.1119/1.10023}
  {\bibfield  {journal} {\bibinfo  {journal} {Am. J. Phys}\ }\textbf {\bibinfo
  {volume} {43}},\ \bibinfo {pages} {22--24} (\bibinfo {year}
  {1975})}\BibitemShut {NoStop}%
\bibitem [{\citenamefont {Salamon}\ and\ \citenamefont
  {Nitzan}(1981)}]{Salamon1981}%
  \BibitemOpen
  \bibfield  {author} {\bibinfo {author} {\bibfnamefont {Peter}\ \bibnamefont
  {Salamon}}\ and\ \bibinfo {author} {\bibfnamefont {Abrahan}\ \bibnamefont
  {Nitzan}},\ }\bibfield  {title} {\enquote {\bibinfo {title} {Finite time
  optimizations of a newton's law carnot cycle},}\ }\href {\doibase
  10.1063/1.441482} {\bibfield  {journal} {\bibinfo  {journal} {J. Chem.
  Phys.}\ }\textbf {\bibinfo {volume} {74}},\ \bibinfo {pages} {3546--3560}
  (\bibinfo {year} {1981})}\BibitemShut {NoStop}%
\bibitem [{\citenamefont {Hille}(2001)}]{Hille2001}%
  \BibitemOpen
  \bibfield  {author} {\bibinfo {author} {\bibfnamefont {Bertil}\ \bibnamefont
  {Hille}},\ }\href
  {https://www.amazon.com/Channels-Excitable-Membranes-Bertil-Hille/dp/0878933212?SubscriptionId=AKIAIOBINVZYXZQZ2U3A&tag=chimbori05-20&linkCode=xm2&camp=2025&creative=165953&creativeASIN=0878933212}
  {\emph {\bibinfo {title} {Ion Channels of Excitable Membranes}}}\ (\bibinfo
  {publisher} {Sinauer Associates is an imprint of Oxford University Press},\
  \bibinfo {year} {2001})\BibitemShut {NoStop}%
\bibitem [{\citenamefont {Berger}(2007)}]{Berger2007}%
  \BibitemOpen
  \bibfield  {author} {\bibinfo {author} {\bibfnamefont {Marcel}\ \bibnamefont
  {Berger}},\ }\href
  {https://www.ebook.de/de/product/2047791/marcel_berger_a_panoramic_view_of_riemannian_geometry.html}
  {\emph {\bibinfo {title} {A Panoramic View of Riemannian Geometry}}}\
  (\bibinfo  {publisher} {Springer Berlin Heidelberg},\ \bibinfo {year}
  {2007})\BibitemShut {NoStop}%
\end{thebibliography}%


\begin{thebibliography}{1}%
\makeatletter
\providecommand \@ifxundefined [1]{%
 \@ifx{#1\undefined}
}%
\providecommand \@ifnum [1]{%
 \ifnum #1\expandafter \@firstoftwo
 \else \expandafter \@secondoftwo
 \fi
}%
\providecommand \@ifx [1]{%
 \ifx #1\expandafter \@firstoftwo
 \else \expandafter \@secondoftwo
 \fi
}%
\providecommand \natexlab [1]{#1}%
\providecommand \enquote  [1]{``#1''}%
\providecommand \bibnamefont  [1]{#1}%
\providecommand \bibfnamefont [1]{#1}%
\providecommand \citenamefont [1]{#1}%
\providecommand \href@noop [0]{\@secondoftwo}%
\providecommand \href [0]{\begingroup \@sanitize@url \@href}%
\providecommand \@href[1]{\@@startlink{#1}\@@href}%
\providecommand \@@href[1]{\endgroup#1\@@endlink}%
\providecommand \@sanitize@url [0]{\catcode `\\12\catcode `\$12\catcode
  `\&12\catcode `\#12\catcode `\^12\catcode `\_12\catcode `\%12\relax}%
\providecommand \@@startlink[1]{}%
\providecommand \@@endlink[0]{}%
\providecommand \url  [0]{\begingroup\@sanitize@url \@url }%
\providecommand \@url [1]{\endgroup\@href {#1}{\urlprefix }}%
\providecommand \urlprefix  [0]{URL }%
\providecommand \Eprint [0]{\href }%
\providecommand \doibase [0]{http://dx.doi.org/}%
\providecommand \selectlanguage [0]{\@gobble}%
\providecommand \bibinfo  [0]{\@secondoftwo}%
\providecommand \bibfield  [0]{\@secondoftwo}%
\providecommand \translation [1]{[#1]}%
\providecommand \BibitemOpen [0]{}%
\providecommand \bibitemStop [0]{}%
\providecommand \bibitemNoStop [0]{.\EOS\space}%
\providecommand \EOS [0]{\spacefactor3000\relax}%
\providecommand \BibitemShut  [1]{\csname bibitem#1\endcsname}%
\let\auto@bib@innerbib\@empty
\bibitem [{\citenamefont {Berger}(2007)}]{Berger2007}%
  \BibitemOpen
  \bibfield  {author} {\bibinfo {author} {\bibfnamefont {M.}~\bibnamefont
  {Berger}},\ }\href
  {https://www.ebook.de/de/product/2047791/marcel_berger_a_panoramic_view_of_riemannian_geometry.html}
  {\emph {\bibinfo {title} {A Panoramic View of Riemannian Geometry}}}\
  (\bibinfo  {publisher} {Springer Berlin Heidelberg},\ \bibinfo {year}
  {2007})\BibitemShut {NoStop}%
\end{thebibliography}%

\end{document}


\title{Supplementary Material: Geodesic bound of the minimum energy expense
to achieve membrane separation within finite time}
\author{Jin-Fu Chen}
\address{Beijing Computational Science Research Center, Beijing 100193, China}
\address{Graduate School of China Academy of Engineering Physics, No. 10 Xibeiwang
East Road, Haidian District, Beijing, 100193, China}
\address{now at: School of Physics, Peking University, Beijing 100871, China}
\author{Ruo-Xun Zhai}
\address{Graduate School of China Academy of Engineering Physics, No. 10 Xibeiwang
East Road, Haidian District, Beijing, 100193, China}
\author{Chang-Pu Sun}
\address{Beijing Computational Science Research Center, Beijing 100193, China}
\address{Graduate School of China Academy of Engineering Physics, No. 10 Xibeiwang
East Road, Haidian District, Beijing, 100193, China}
\author{Hui Dong}
\email{hdong@gscaep.ac.cn}

\address{Graduate School of China Academy of Engineering Physics, No. 10 Xibeiwang
East Road, Haidian District, Beijing, 100193, China}
\date{\today}

\maketitle
The document is devoted to providing the detailed derivation and the
discussions to the main text.

\section{First-order approximation for evolution equation}

The\textbf{ }quadratic form of $\dot{W}_{\mathrm{ex}}$ is obtained
by keeping the series solutions of evolution equations to the first-order
term of $1/\tau$. The approximation is validated when the time of
process $\tau$ is much larger than the relaxation time of the system
$\tau_{\alpha,\beta}$ and $\tau_{h}$, or in the endo-reversible
region. For a given protocol $\tilde{x}_{k}(s),$ $V_{k}(t)$ is obtained
by $V_{k}(t)=V_{t}\tilde{x}_{k}(t/\tau)$ , ($k=L,R$). The time evolution
of the system are described by Eqs. (1) - (3). We define the variation
of $N_{\sigma k}$ and $T$ from the corresponding equilibrium states
as

\begin{equation}
\begin{split}\Delta N_{\alpha L}\equiv & N_{\alpha L}-\frac{V_{L}}{V_{M}+V_{L}}N_{\alpha}\,,\quad\Delta N_{\beta R}\equiv N_{\beta R}-\frac{V_{R}}{V_{M}+V_{R}}N_{\beta}\,,\quad\Delta T\equiv T-T_{0}\,.\end{split}
\label{eq:def_Deltas}
\end{equation}
The evolution equations of the variations are
\begin{align}
\frac{d}{dt}(\Delta N_{\alpha L})= & -N_{\alpha}\frac{d}{dt}\left(\text{\ensuremath{\frac{V_{L}}{V_{M}+V_{L}}}}\right)-\mu_{\alpha A}\mathcal{A}\left(\frac{1}{V_{M}}+\frac{1}{V_{L}}\right)\Delta N_{\alpha L}\,,\label{eq:Del_N_alphaL}\\
\frac{d}{dt}(\Delta N_{\beta R})= & -N_{\beta}\frac{d}{dt}\left(\frac{V_{R}}{V_{M}+V_{R}}\right)-\mu_{\beta B}\mathcal{A}\left(\frac{1}{V_{M}}+\frac{1}{V_{R}}\right)\Delta N_{\beta R}\,,\label{eq:Del_N_betaR}\\
\frac{d}{dt}(\Delta T)= & -\gamma\Delta T+\frac{k_{B}(T_{0}+\Delta T)}{C_{V}}f(t)\,,\label{eq:Del_T}
\end{align}
where $f(t)\equiv\dot{V}_{L}\Big[N_{\beta}/(V_{M}+V_{R})-\left(1/V_{M}+1/V_{L}\right)\Delta N_{\alpha L}-\Delta N_{\beta R}/V_{M}\Big]+\dot{V}_{R}\Big[N_{\alpha}/(V_{M}+V_{L})-\left(1/V_{M}+1/V_{R}\right)\Delta N_{\beta R}-\Delta N_{\alpha L}/V_{M}\Big]\,.$
To the first-order term of $1/\tau\,$, Eqs. (\ref{eq:Del_N_alphaL})-(\ref{eq:Del_T})
are solved by self-iteration method,
\begin{equation}
\begin{split}\Delta N_{\alpha L}(t)\approx & -\frac{N_{\alpha}}{\mu_{\alpha A}\mathcal{A}\left(\frac{1}{V_{M}}+\frac{1}{V_{L}}\right)}\frac{d}{dt}\left(\frac{V_{L}}{V_{t}-V_{R}}\right)\,,\\
\Delta N_{\beta R}(t)\approx & -\frac{N_{\beta}}{\mu_{\beta B}\mathcal{A}\left(\frac{1}{V_{M}}+\frac{1}{V_{R}}\right)}\frac{d}{dt}\left(\frac{V_{R}}{V_{t}-V_{L}}\right)\,,\\
\Delta T(t)\approx & \frac{k_{B}T_{0}}{\gamma C_{V}}\left(\frac{\dot{V}_{L}}{V_{t}-V_{L}}N_{\beta}+\frac{\dot{V}_{R}}{V_{t}-V_{R}}N_{\alpha}\right)\,.
\end{split}
\label{eq:firstorder_solution}
\end{equation}
The details of iteration method are presented in the Supplementary
Information.

\section{Excess work rate}

The work rate of the system is $\dot{W}=\sum_{\sigma=\alpha,\beta}\sum_{k=L,M,R}p_{\sigma k}\dot{V}_{k}=-p_{\alpha L}\dot{V}_{L}-p_{\beta R}\dot{V}_{R}-(p_{\alpha M}+p_{\beta M})(-\dot{V}_{L}-\dot{V}_{R})\,$.
Together with Eq. (1), the work rate is expressed as
\begin{align}
\dot{W}= & k_{B}T\dot{V}_{L}\left(\frac{N_{\alpha M}}{V_{M}}+\frac{N_{\beta M}}{V_{M}}-\frac{N_{\alpha L}}{V_{L}}\right)+k_{B}T\dot{V}_{R}\left(\frac{N_{\alpha M}}{V_{M}}+\frac{N_{\beta M}}{V_{M}}-\frac{N_{\beta R}}{V_{R}}\right)\,.\label{eq:Power_tot}
\end{align}
Substituting Eq. (\ref{eq:firstorder_solution}) into Eq. (\ref{eq:Power_tot})
and keeping to quadratic term of $1/\tau$, we obtain the leading
term of the total work rate, which is divided into the reversible
part and excess part, i.e., $\dot{W}=\dot{W}^{(0)}+\dot{W}_{\mathrm{ex}}\,$.
In reversible separation processes, the system stays at a transient
equilibrium state at every moment, the partial pressures of type-$\alpha$
(type-$\beta$) molecules over the A (B) membrane are equal $p_{\alpha L}=p_{\alpha M}\,,p_{\beta R}=p_{\beta M}$,
then all the work done during this process is reversible, and the
work rate is
\[
\dot{W}^{(0)}=k_{B}T_{0}\left(\dot{V}_{L}\cdot\frac{N_{\beta}}{V_{t}-V_{L}}+\dot{V}_{R}\cdot\frac{N_{\alpha}}{V_{t}-V_{R}}\right)\,.
\]
The excess work rate is obtained by subtracting the reversible one
$\dot{W}^{(0)}$ from $\dot{W}\,$, i.e., $\dot{W}_{\mathrm{ex}}=\dot{W}-\dot{W}^{(0)}$,
and the result is given by Eq. (5).

\section{Minimal reversible work}

For the quasi-static isothermal separation with infinite operation
time, the partial pressure of molecules in each compartment satisfies
$p_{\alpha M}=p_{\alpha L}$ and $p_{\beta M}=p_{\beta R}$. The molecular
number of the gas in each compartment is
\begin{align}
N_{\alpha L} & =N_{\alpha}\frac{V_{L}}{V_{L}+V_{M}}\,,\label{eq:N110}\\
N_{\beta R} & =N_{\beta}\frac{V_{R}}{V_{R}+V_{M}}\,,\label{eq:N220}
\end{align}
with $N_{\alpha M}=N_{\alpha}-N_{\alpha L}$ and $N_{\beta M}=N_{\beta}-N_{\beta R}$
. The temperature of the system remains the same as that of the bath
$T=T_{0}$. The partial pressure $p_{\alpha M}=N_{\alpha}k_{B}T_{0}/(V_{M}+V_{L})$
and $p_{\beta M}=N_{\beta}k_{B}T_{0}/(V_{M}+V_{R})$ in the quasi-static
process only depends on the volumes $V_{k}$ ($k=L,M,R$) of the compartments.
We set the final volumes of the left and the right compartments to
be $V_{L}(\tau)$ and $V_{R}(\tau)$ with $V_{L}(\tau)+V_{R}(\tau)=V_{t}$.
The performed work for the quasi-static process is

\begin{align}
W^{(0)} & =-k_{B}T_{0}(N_{\alpha}\ln\frac{V_{L}(\tau)}{V_{t}}+N_{\beta}\ln\frac{V_{R}(\tau)}{V_{t}}).
\end{align}
The minimal quasi-static performed work is $W_{(0)}^{\mathrm{min}}=-N_{t}k_{B}T_{0}[\epsilon_{\alpha}\ln\epsilon_{\alpha}+\epsilon_{\beta}\ln\epsilon_{\beta}]$,
which is reached by choosing $V_{L}(\tau)/V_{t}=\epsilon_{\alpha}\equiv N_{\alpha}/N$
and $V_{R}(\tau)/V_{t}=\epsilon_{\beta}\equiv N_{\beta}/N$. For the
situation in the main text where $N_{t}=2N_{A},\epsilon_{\alpha}=\epsilon_{\beta}=0.5,T_{0}=298.15\mathrm{K}$,
the minimal reversible work is $W^{(0)}=3.44\times10^{3}\mathrm{J}$.

\section{Self-iteration method to solve differential equations}

We employ the self-iteration method to Eqs. (15)-(17) of the main
text. We rewrite Eq. (15) into
\[
\Delta N_{\alpha L}=-\frac{N_{\alpha}}{\mu_{\alpha}\mathcal{A}\left(\frac{1}{V_{M}}+\frac{1}{V_{L}}\right)}\frac{d}{dt}\left(\frac{V_{L}}{V_{M}+V_{L}}\right)-\frac{1}{\mu_{\alpha}\mathcal{A}\left(\frac{1}{V_{M}}+\frac{1}{V_{L}}\right)}\frac{d}{dt}\left(\Delta N_{\alpha L}\right)\,.
\]
The left hand side is the variable itself, while its time derivative
is on the right hand side. By self-iterating the equation, we obtain
a solution in the series form
\[
\begin{split}\Delta N_{\alpha L}= & -\frac{N_{\alpha}}{\mu_{\alpha}\mathcal{A}\left(\frac{1}{V_{M}}+\frac{1}{V_{L}}\right)}\frac{d}{dt}\left(\frac{V_{L}}{V_{M}+V_{L}}\right)\\
 & -\frac{N_{\alpha}}{\mu_{\alpha}\mathcal{A}\left(\frac{1}{V_{M}}+\frac{1}{V_{L}}\right)}\frac{d}{dt}\left\{ -\frac{N_{\alpha}}{\mu_{\alpha}\mathcal{A}\left(\frac{1}{V_{M}}+\frac{1}{V_{L}}\right)}\frac{d}{dt}\left(\frac{V_{L}}{V_{M}+V_{L}}\right)-\frac{1}{\mu_{\alpha}\mathcal{A}\left(\frac{1}{V_{M}}+\frac{1}{V_{L}}\right)}\frac{d}{dt}\left(\Delta N_{\alpha L}\right)\right\} \\
= & N_{\alpha}\sum_{n=1}^{\infty}\left[-\frac{1}{\mu_{\alpha}\mathcal{A}\left(\frac{1}{V_{M}}+\frac{1}{V_{L}}\right)}\frac{d}{dt}\right]^{n}\left(\frac{V_{L}}{V_{M}+V_{L}}\right)\,.
\end{split}
\]
We choose a given protocol to control the volume $V_{k}(t)=V_{t}\tilde{x}_{k}(t/\tau)$
with the operation time $\tau$ of the whole separation process. The
time derivative is substituted into 
\begin{equation}
\frac{d^{n}}{dt^{n}}x_{k}(t)=\frac{d^{n}}{dt^{n}}\tilde{x}_{k}(t/\tau)=\frac{1}{\tau^{n}}\frac{d^{n}}{ds^{n}}\tilde{x}_{k}(s)\,,\label{eq:9}
\end{equation}
with the rescaled time $s=t/\tau$. Therefore, each time derivative
contributes a $1/\tau$ factor. For slow separation process, we keep
to the lowest-order term of $1/\tau$, and obtain an approximate solution
in Eq. (18) of the main text. The approximate solutions to $\Delta N_{\beta R}$
and $\Delta T$ can be obtained similarly.

\section{Solving the geodesic equations}

In this section, we show the method to numerically evaluate the geodesics
paths connecting $(0,0)$ and $(\epsilon_{\alpha},\epsilon_{\beta})$.
We need to solve the geodesic equations
\begin{equation}
\frac{d^{2}x^{k}}{dr^{2}}+\Gamma_{\;ij}^{k}\frac{dx^{i}}{dr}\frac{dx^{j}}{dr}=0\,,\label{eq:Geodesic}
\end{equation}
where $(x^{1},x^{2})$ represents the coordinate $(x_{L},x_{R})$
in the Riemann configuration space. The repeated indexes are summed.
The expressions of elements of the metric are explicitly
\begin{equation}
\begin{split}g_{11}= & N_{t}k_{B}T_{0}\left[\epsilon_{\beta}\tau_{\beta}(x^{2})^{2}/(1-x_{1})^{3}+\epsilon_{\alpha}\tau_{\alpha}/(1-x_{2})+2/3\cdot\epsilon_{\beta}^{2}\tau_{h}/(1-x_{1})^{2}\right]\,,\\
g_{12}=g_{21}= & N_{t}k_{B}T_{0}\left[\epsilon_{\alpha}\tau_{\alpha}x^{1}/(1-x^{2})^{2}+\epsilon_{\beta}\tau_{\beta}x^{2}/(1-x^{1})^{2}+2/3\cdot\epsilon_{\alpha}\epsilon_{\beta}\tau_{h}/(1-x^{1})(1-x^{2})\right]\,,\\
g_{22}= & N_{t}k_{B}T_{0}\left[\epsilon_{\alpha}\tau_{\alpha}(x^{1})^{2}/(1-x^{2})^{3}+\epsilon_{\beta}\tau_{\beta}/(1-x^{1})+2/3\cdot\epsilon_{\alpha}^{2}\tau_{h}/(1-x^{2})^{2}\right]\,.
\end{split}
\end{equation}
The Christoffel symbols are obtained as
\begin{align}
\Gamma_{111}= & \frac{N_{t}k_{B}T_{0}}{2}\left[\frac{4}{3}\cdot\frac{\epsilon_{\beta}^{2}\tau_{h}}{(1-x^{1})^{3}}+\frac{3\epsilon_{\beta}\tau_{\beta}(x^{2})^{2}}{(1-x^{1})^{4}}\right]\,,\\
\Gamma_{222}= & \frac{N_{t}k_{B}T_{0}}{2}\left[\frac{4}{3}\cdot\frac{\epsilon_{\alpha}^{2}\tau_{h}}{(1-x^{2})^{3}}+\frac{3\epsilon_{\alpha}\tau_{\alpha}(x^{1})^{2}}{(1-x^{2})^{4}}\right]\,,\\
\Gamma_{121}=\Gamma_{211}= & \frac{N_{t}k_{B}T_{0}}{2}\left[\frac{\epsilon_{\beta}\tau_{\beta}}{(1-x^{1})^{2}}+\frac{2\epsilon_{\alpha}\tau_{\alpha}x^{1}}{(1-x^{2})^{3}}\right]\,,\\
\Gamma_{122}=\Gamma_{212}= & \frac{N_{t}k_{B}T_{0}}{2}\left[\frac{\epsilon_{\alpha}\tau_{\alpha}}{(1-x^{2})^{2}}+\frac{2\epsilon_{\beta}\tau_{\beta}x^{2}}{(1-x^{1})^{3}}\right]\,,\\
\Gamma_{221}= & \frac{N_{t}k_{B}T_{0}}{2}\left[\frac{\epsilon_{\beta}\tau_{\beta}}{(1-x^{1})^{2}}+\frac{2\epsilon_{\alpha}\tau_{\alpha}x^{1}}{(1-x^{2})^{3}}+\frac{4}{3}\cdot\frac{\epsilon_{\alpha}\epsilon_{\beta}\tau_{h}}{(1-x^{1})(1-x^{2})^{2}}\right]\,,\\
\Gamma_{112}= & \frac{N_{t}k_{B}T_{0}}{2}\left[\frac{\epsilon_{\alpha}\tau_{\beta}}{(1-x^{2})^{2}}+\frac{2\epsilon_{\beta}\tau_{\beta}x^{2}}{(1-x^{1})^{3}}+\frac{4}{3}\cdot\frac{\epsilon_{\alpha}\epsilon_{\beta}\tau_{h}}{(1-x^{2})(1-x^{1})^{2}}\right]\,.
\end{align}
In Eq. (\ref{eq:Geodesic}), we have $\Gamma_{\ jk}^{i}\equiv g^{li}\Gamma_{jkl}\,$
with the inverse metric $g^{ij}$. We apply the shooting method \citep{Berger2007}
in searching of the geodesic path with a target point $(\epsilon_{\alpha},\epsilon_{\beta})$.
The boundary value problem of the ordinary differential equation is
converted into an initial value problem by adding the initial derivative
$x^{i\prime}(0)$ into the initial condition. A proper value of $x^{i\prime}(0)$
is searched to allow the end of solved geodesics path match the target
point $(\epsilon_{\alpha},\epsilon_{\beta})$. We represent the solved
geodesics path as a function of the arc length $r$, i.e., $x^{i}=x^{i}(r),\,0<r<\mathcal{L}\,$,
and $\mathcal{L}$ is the thermodynamics length of the path. To represent
the protocol $\tilde{x}_{L,R}(s)\,$, we can re-parameterize $x^{i}(r)$
by $s\equiv r/\mathcal{L}$, and the protocol is explicitly $\tilde{x}_{L}(s)=x^{1}(s\mathcal{L}),\tilde{x}_{R}(s)=x^{2}(s\mathcal{L})\,$.

It is worth noting that the metric at $(\epsilon_{\alpha},\epsilon_{\beta})$
is degenerate
\[
G=N_{t}k_{B}T_{0}\left(\tau_{\alpha}+\tau_{\beta}+\frac{2}{3}\cdot\tau_{h}\right)\begin{pmatrix}1 & 1\\
1 & 1
\end{pmatrix}\,,
\]
inducing the divergence of $G^{-1}\,$. And the Christoffel symbols
$\Gamma_{\ ij}^{k}$ are also divergent at the target point $(\epsilon_{\alpha},\epsilon_{\beta})$.
When numerically solving the geodesic equations, we set the target
point $(\epsilon_{\alpha},\epsilon_{\beta})=(0.5,0.5)$. Due to the
divergence of $\Gamma_{\mu\nu}^{\rho}$, the integral step length
$\delta r$ is supposed to decrease to meet the accuracy requirement.
The solving procedure stops once the required step length is smaller
than a minimum one, set by our algorithm. There is a small distance
between the stop point and the target point. In the numerical calculation
of the geodesic equation, we set the minimum step length to be $\delta r_{\mathrm{min}}=10^{-13}(\mathrm{J\cdot s})^{1/2}$,
and obtain a stop point with a Euclidean distance $9.7\times10^{-8}$
away from $(0.5,0.5)$, or a geometric distance $\Delta\mathcal{L}\approx4.9\times10^{-6}(\mathrm{J\cdot s})^{1/2}$
(the relative error of the thermodynamic length of the whole path
is $\Delta\mathcal{L}/\mathcal{L}_{\mathrm{geo}}\approx8.2\times10^{-8}$).
This distance is small enough to be ignored.

\section{Singularities of the evolution equations}

With the protocol $\tilde{x}_{L,R}(s)$ obtained from the geodesic
path, we can numerically simulate the evolution of the ideal gas by
the evolution Eqs. (15)-(17) of the main text. With a certain operation
time $\tau$, the volumes are expressed explicitly $V_{k}(t)=V_{t}\tilde{x}_{k}(t/\tau)$
as functions of time. The evolution equations have two removable singularities
at the beginning $t=0$ and the end $t=\tau$ of the separation process,
since the particle densities $N_{\sigma k}/V_{k}$ become in the form
$0/0$. The singularities will bring computational difficulty for
the numerical simulation. Here we illustrate how to avoid the singularities
in the numerical simulation without affecting the final result.

At the beginning $t=0$ of the separation process, the particle densities
$N_{\alpha L}/V_{L}$ and $N_{\beta R}/V_{R}$ of the left and the
right compartments lead to a removable singularity of the evolution
equations since both the numerators and denominators approach $0$.
We consider an initial equilibrium state of the gas, and thus the
initial values of the molecular numbers are given by Eqs. (\ref{eq:N110})
and (\ref{eq:N220}). In the numerical simulation, this singularity
can be avoided by choosing an initial time slightly larger than zero,
such as $\delta t_{i}=\tau\times10^{-8}$. The initial $V_{k}$ is
slightly larger than zero by $V_{k}(\delta t_{i})$, and the initial
$N_{\sigma k}$ is set according to Eqs. (\ref{eq:N110}) and (\ref{eq:N220}).
We also remark that in real situation the volumes cannot be perfectly
zero.

At the end $t=\tau$ of the separation process, the particle densities
$(N_{\alpha}-N_{\alpha L})/(V_{t}-V_{L}-V_{R})$ and $(N_{\beta}-N_{\beta R})/(V_{t}-V_{L}-V_{R})$
of the middle compartment lead to a removable singularity of the evolution
equations. The final point obtained by the shooting method is not
exactly $(0.5,0.5)$. And the volume of the middle compartment has
a small positive value at the end of the separation process. Thus,
the differential equations can be solved numerically with the protocol
obtained from the shooting method. 

\section{Excess work for short operation time}

\begin{figure}
\begin{centering}
\includegraphics{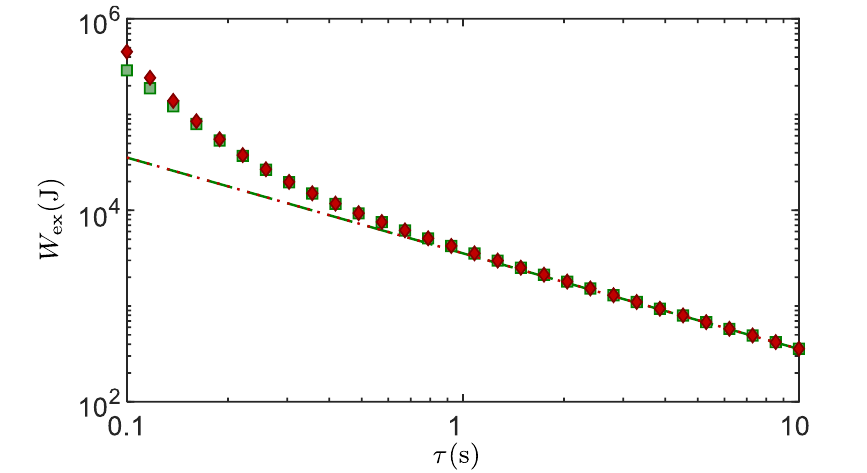}
\par\end{centering}
\caption{Excess work $W_{\mathrm{ex}}$ versus operation time $\tau\,$ in
the short time region. The red diamonds (green squares) represent
the irreversible work obtained by solving the evolution equation of
the system along the symmetric (symmetry-breaking) path. The red dashed
line (green solid line) represents the excess work predicted by the
thermodynamic length $W_{\mathrm{ex}}=\mathcal{L}^{2}/\tau$. \label{fig:shottimeexcesswork}}
\end{figure}

The above discussion on the excess work is for the long time region
where the operation time $\tau$ is much larger than the time scale
of relaxation. In the short-time region, we show the excess work in
Fig. \ref{fig:shottimeexcesswork} with numerical simulation. The
red diamonds (green squares) represent the irreversible work obtained
by solving the evolution equation of the system for the symmetric
(symmetry-breaking) protocol. The red dashed line (green solid line)
represents the excess work predicted by the thermodynamic length $W_{\mathrm{ex}}=\mathcal{L}^{2}/\tau$.
The excess work is far larger than the reversible work $W^{(0)}$
in Sec. II.

\bibliographystyle{apsrev4-1}
\bibliography{separationbib}